\documentclass[12pt]{article}
\usepackage{mathptmx}
\usepackage{graphicx}
\usepackage{amsmath}
\usepackage{amsfonts}
\newtheorem{thm}{Theorem}[section]
\newtheorem{lem}{Lemma}[section]

\newtheorem{Corollary}{Corollary}[section]

\newtheorem{tb}{Table}

\title{ Direct and inverse results for Kantorovich type exponential sampling series}
\author{Sathish Kumar Angamuthu \thanks{Department of Mathematics, Visvesvaraya National Institute of Technology, Nagpur, Nagpur-440010, India. \newline E-mail: mathsatish9@gmail.com}
 \and
  Shivam Bajpeyi \thanks{Department of Mathematics, Visvesvaraya National Institute of Technology, Nagpur, Nagpur-440010, India. \newline E-mail: shivambajpai1010@gmail.com}
  }
\date{}

\begin{document}
\maketitle
\bibliographystyle{plain}
\abstract{In this article, we analyze the behaviour of the new family of Kantorovich type exponential sampling series. We obtain the point-wise approximation theorem and Voronovskaya type theorem for the series $(I_{w}^{\chi})_{w>0}.$ Further, we obtain a representation formula and an inverse result approximation for these operators. Finally, we give some examples of kernel functions to which the theory can be applied along with the graphical representation.

\endabstract



\section{Introduction}
The exponential sampling methods play the key role in solving problems in the area of optical physics and engineering, precisely in the phenomena like Fraunhofer diffraction, light scattering etc \cite{casasent,bertero,gori,ostrowsky}. It all began when a group of optical physicists and engineers Bartero, Pike \cite{bertero} and Gori \cite{gori} presented a representation formula known as \textit{exponential sampling formula}, for the class of Mellin band-limited function having exponentially spaced sample points, which is also considered as the Mellin-version of the well known \textit{Shannon sampling theorem} ( see \cite{butzer2}). But, the pioneering idea of mathematical study of exponential sampling formula is credited to Butzer and Jansche. Butzer et al. \cite{butzer5} proved the exponential sampling formula mathematically using the theory of Mellin transform and Mellin approximation, which was first studied separately by Mamedov \cite{mamedeo} and then developed by Butzer and Jansche \cite{butzer3,butzer4,butzer5,butzer7}. We mention some of the work related to the theory of Mellin transform and Mellin approximation (see \cite{bardaro1,bardaro2,bardaro3}). Bardaro et al.\cite{bardaro7} made a significant development in this direction when they replaced the $lin_{c}$ function in the exponential sampling formula by more general kernel function. Let $ x \in \mathbb{R}^{+}$ and $ w >0.$ Then, the generalized exponential sampling series (see \cite{bardaro7}) is defined by
\begin{equation} \label{classical}
(S_{w}^{\chi}f)(x)= \sum_{k=- \infty}^{+\infty} \chi(e^{-k} x^{w}) f( e^{\frac{k}{w}})
\end{equation}
where $ f:\mathbb{R}^{+} \rightarrow \mathbb{R} $ be any function for which the series is absolutely convergent. The convergence of the above series studied in Mellin-Lebesgue spaces in \cite{bardaro11}. The importance of the above series appears when we need to reconstruct a not necessarily Mellin band-limited signal in approximate sense having samples which are exponentially spaced. This provides us a useful tool to approximate a signal by using its values at the node $ (e^{\frac{k}{w}}).$ But practically it is difficult to have the exact sample value at the node $ (e^{\frac{k}{w}})$ always. To overcome this problem, one can replace the value $ f( e^{\frac{k}{w}})$ by the mean value of $ f(e^{x})$ in the interval $ \big[ \frac{k}{w}, \frac{k+1}{w} \big],$ for $ k \in \mathbb{Z}, w>0.$ This idea motivates us to define the Kantorovich-version of the generalized exponential sampling series (\ref{classical}) and this work is inspired from the classical Kantorovich version (\cite{bardaro10,bardaro5,bardaro6,costa2,costa3,costa1,vinti2,vinti1,vinti3,orlova}) of the generalized sampling series introduced by Butzer in \cite{butzer2} .

The study of Kantorovich type generalizations of approximation operators is an important subject in approximation theory, as they can be used to approximate
Lebesgue integrable functions. In the last few decades, Kantorovich modifications of several operators have been constructed and their approximation behavior studied,
we mention some of the work in this direction e.g., \cite{gupta1,maria,PNA,Agrawal1,Agrawal2,Tuncer,Tuncer1,kajla,vijay1,vijay2} etc.

Let $C(\mathbb{R}^+$) be the space of all continuous and bounded functions on $\mathbb{R}^+$.
A function $f \in C(\mathbb{R}^+$) is called log-uniformly continuous on $\mathbb{R}^+$, if for any given
$\epsilon > 0,$ there exists $\delta > 0$ such that $|f(u) -f(v)| < \epsilon$ whenever $| \log u - \log v| \leq \delta,$ for any $u, v \in \mathbb{R}^{+}.$ We denote the space of all   log-uniformly continuous functions defined on  $\mathbb{R}^{+}$ by $\mathcal{C}(\mathbb{R}^+).$ We consider $M(\mathbb{R}^{+})$ as the class of all Lebesgue measurable functions on $\mathbb{R}^+$ and $L^{\infty}(\mathbb{R}^{+}) $ as the space of all bounded functions on  $\mathbb{R}^{+}$ throughout this paper.

For $1 \leq p < +\infty$, let $ L^p(\mathbb{R}^+$) be the space of all the Lebesgue measurable and $p$-integrable functions defined on $\mathbb{R}^+$ equipped with the usual norm $\Vert f \Vert_p$. For $c \in \mathbb{R}$, we define the space
$$X_c = \{f : \mathbb{R}^+ \rightarrow \mathbb{C} : f(\cdot)(\cdot)^{c-1} \in  L^1(\mathbb{R}^+)\}$$ equipped with the norm
$$\Vert f \Vert_{X_c} = \Vert f(\cdot)(\cdot)^{c-1} \Vert_1 = \int_0^{+\infty} |f(u)|u^{c-1}du.  $$

The Mellin transform of a function $f \in X_c$ is defined by
$$\hat{M}[f](s) := \int_0^{+\infty} u^{s-1}f(u)\ du \ , \,  \ (s = c + it, t \in \mathbb{R}).$$
A function $f \in X_{c} \cap C(\mathbb{R}^+), c \in  \mathbb{R}$ is called Mellin
band-limited in the interval $[-\eta, \eta],$ if $\hat{M}[f](c+iw) = 0$ for all $|w| > \eta ,\ \eta \in \mathbb{R}^+.$

Let $f : \mathbb{R}^+ \rightarrow \mathbb{C}$ and $c \in \mathbb{R}.$ Then, Mellin differential operator $\theta_c$ is defined by
$$\theta_cf(x) := xf'(x) + cf(x), \ \ \ \ \ x \in  \mathbb{R}^+ .$$
We consider $ \theta f(x) := \theta_{0} f(x) $ throughout this paper.
The Mellin differential operator of order $r \in \mathbb{N}$ is
defined by $\theta_c^1 := \theta_c, \, \, \  \theta_c^r = \theta_c(\theta_c^{r-1}).$

Let $ \chi :\mathbb{R}^{+} \rightarrow \mathbb{R}$ be the kernel function which is continuous on $\mathbb{R}^{+}$ such that it satisfies the following conditions: \\
(i) For every $ x \in \mathbb{R}^{+},$
$$  \sum_{k=- \infty}^{+\infty} \chi(e^{-k} x^{w}) =1.$$\\
(ii) $M_{2}(\chi) <+ \infty$  and
$$ \lim_{\gamma \rightarrow + \infty} \sum_{|k-\log(u)|>  \gamma} |\chi(e^{-k} u)| \ |k- \log(u)|^{2}=0,$$ uniformly with respect to $ u \in \mathbb{R}^{+}.$\\
We define the algebraic moments of order $\nu$ for the kernel function $ \chi$ as
$$ m_{\nu}(\chi,u):= \sum_{k= - \infty}^{+\infty}  \chi(e^{-k} u) (k- \log(u))^{\nu}, \hspace{0.5cm} \forall \ u \in \mathbb{R}^{+}.$$
Similarly, the absolute moment of order $\nu$ can be defined as
$$ M_{\nu}(\chi,u):= \sum_{k= - \infty}^{+\infty}  |\chi(e^{-k} u)| |k- \log(u)|^{\nu},  \hspace{0.5cm} \forall \ u \in \mathbb{R}^{+}.$$

We define $ \displaystyle M_{\nu}(\chi):= \sup_{u \in \mathbb{R}^{+}} M_{\nu}(\chi,u). $

Let $ x \in \mathbb{R}^{+}$ and $ w >0.$ We define the Kantorovich version of the exponential sampling series given by (\ref{classical}) as follows:
\begin{equation} \label{main}
(I_{w}^{\chi}f)(x)= \sum_{k= - \infty}^{+\infty} \chi(e^{-k} x^{w})\  w \int_{\frac{k}{w}}^{\frac{k+1}{w}} f(e^{u})\  du, \ \
\end{equation}
where $ f: \mathbb{R}^{+} \rightarrow \mathbb{R}$ is locally integrable such that the above series is convergent for every $ x \in \mathbb{R}^{+}$. It is clear that for $f \in L^{\infty}(\mathbb{R}^{+}),$ the above series is well defined for every $x \in \mathbb{R}^{+}.$ The paper is organized as follows. In section 2 we prove the point-wise approximation theorem and Voronovskaya type asymptotic formula for the Kantorovich exponential sampling operators. In section 3 we obtain the representation formula which gives the relationship between the exponential sampling operators with Kantorovich exponential sampling operators. The main result of this section is inverse result for the operators (\ref{main}). Finally we verify the assumptions used in the theory by using Mellin $B$-splines and Mellin's Fejer kernel and show the approximation of functions by (\ref{main}) graphically.

\section{Approximation Results}

In this section, we obtain some direct results e.g. pointwise convergence theorem and Voronovskaya type asymptotic formula for the Kantorovich exponential sampling operators (\ref{main}).

\begin{thm}\label{theorem1}
Let $ f\in M(\mathbb{R}^{+}) \cap L^{\infty}(\mathbb{R}^{+}).$ Then, the series (\ref{main}) converges to $f(x)$ at every point $x \in \mathbb{R}^{+},$ the point of continuity of $f$. Moreover, for $ f \in C(\mathbb{R}^{+})$
$$ \lim_{w \rightarrow \infty} \|I_{w}^{\chi}f - f \|_{\infty} = 0.$$
\end{thm}

\noindent\bf{Proof.}\rm \
Using the condition (i), we obtain
\begin{eqnarray*}
|I_{w}^{\chi}f(x)-f(x)| &=& \bigg| \sum_{k= - \infty}^{+\infty} \chi(e^{-k} x^{w})\  w \int_{\frac{k}{w}}^{\frac{k+1}{w}} (f(e^{u}) - f(x))du \bigg| \\
&\leq& \left( \sum_{\big|k-w \log(x)\big|< \frac{w \delta}{2}} +\sum_{\big|k- w\log(x)\big| \geq \frac{w \delta}{2}} \right) \big| \chi(e^{-k} x^{w})\big| \ w \int_{\frac{k}{w}}^{\frac{k+1}{w}} |f(e^{u}) - f(x)| \ du. \\
&:=& I_{1}+I_{2}. \
\end{eqnarray*}
Since $ f \in C(\mathbb{R}^{+}),$ for every $ \epsilon >0$ there exists $\delta>0$ such that $ |f(e^{u}) - f(x)|< \epsilon,$ whenever $|u- \log(x)|< \delta.$
Let $w^{'}$ be fixed in such a way that $ \frac{1}{w}< \frac{\delta}{2} $ for every $ w > w^{'}.$ Now, for  $ u \in\big[ \frac{k}{w}, \frac{k+1}{w} \big]$ and $ w > w^{'}$ and we have
$$ |u- \log(x)| \leq \Big | u - \frac{k}{w} \Big|+\Big |\frac{k}{w} - \log(x) \Big| \leq \delta,$$
whenever $ \big|\frac{k}{w} - \log(x)\big|< \frac{\delta}{2}.$
This gives $|I_{1}| < \epsilon M_{0}(\chi).$ Similarly, we estimate $ I_{2}.$
\begin{eqnarray*}
|I_{2}| &\leq& \ 2 \|f\|_{\infty} \sum_{\big|k- w \log(x)\big| \geq w \delta} \big| \chi(e^{-k} x^{w})\big|\leq 2 \|f\|_{\infty} \  \epsilon.
\end{eqnarray*}
Combining the estimates of $I_{1}-I_{2},$ we get the  desired result.\\

Next, we derive the following asymptotic formula for the operators $(I_{w}^{\chi}f)_{w>0} .$
\begin{thm}\label{theorem2}
Let $ f \in C^{(2)}(\mathbb{R}^{+})$ and  $\chi$ be the kernel function such that its first order moment vanishes for all $ u \in \mathbb{R}^{+}.$ Then, we have
$$ \lim_{w\rightarrow \infty} w \big[ (I_{w}^{\chi}f)(x) - f(x) \big] = \frac{(\theta f)(x)}{2}.$$
 \end{thm}

\noindent\bf{Proof.}\rm \ For $f \in C^{2}(\mathbb{R}^{+}),$ the Taylor's formula in terms of Mellin derivatives (\cite{butzer3}) upto second order term can be written as
$$f(e^{u})=f(x)+ (\theta f)(x) (u-\log(x)) + \frac{(\theta^{2} f)(x)}{2!} (u-\log(x))^{2} + h \Big(\frac{e^{u}}{x}\Big) (u-\log(x))^{2},$$ where $h$ is a bounded function such that $\displaystyle \lim_{t \rightarrow 1} h(t)=0.$ In view of (\ref{main}) we obtain
\begin{eqnarray*}
[(I_{w}^{\chi}f)(x) - f(x)] &=& \sum_{k= - \infty}^{+\infty} \chi(e^{-k} x^{w})\  w \int_{\frac{k}{w}}^{\frac{k+1}{w}} \big[ (\theta f)(x) (u-\log(x)) + \\&&
\frac{(\theta^{2} f)(x)}{2!} (u-\log(x))^{2} + h \Big(\frac{e^{u}}{x}\Big) (u-\log(x))^2 \Big] \ du \\
&:=& I_{1}+I_{2}+ I_{3}.
\end{eqnarray*}
First we evaluate $I_{1}.$
\begin{eqnarray*}
I_{1}&=& \sum_{k= - \infty}^{+\infty}  \chi(e^{-k} x^{w})\  w \int_{\frac{k}{w}}^{\frac{k+1}{w}} \big[ (\theta f)(x) (u-\log(x))\big]\ du \\
&=& (\theta f)(x) \sum_{k= - \infty}^{+\infty} \chi(e^{-k} x^{w})\ w \int_{\frac{k}{w}}^{\frac{k+1}{w}} (u-\log(x)) \ du=\frac{(\theta f)(x)}{2w} .
\end{eqnarray*}
Now, we estimate $I_{2}.$
\begin{eqnarray*}
I_{2}&=& \sum_{k= - \infty}^{+\infty}  \chi(e^{-k} x^{w})\  w \int_{\frac{k}{w}}^{\frac{k+1}{w}} \Big[ \frac{(\theta^{(2)} f)(x)}{2!} (u-\log(x))^{2} \Big]\ du \\
&=& \frac{(\theta^{(2)} f)(x)}{6} \sum_{k= - \infty}^{+\infty} \chi(e^{-k} x^{w})\ w \int_{\frac{k}{w}}^{\frac{k+1}{w}} (u-\log(x))^{2} \ du\\
&=& \frac{(\theta^{(2)} f)(x)}{6 w^{2}} \big( 1 + 3 m_{1}(\chi,u) + 3 m_{2}(\chi,u) \big).
\end{eqnarray*}
From the above estimates $I_1$ and $I_2,$ we obtain $w(I_1+I_2)\rightarrow \frac{(\theta f)(x)}{2}$ as $w\rightarrow \infty.$
Now $I_{3}$ can be written as
\begin{eqnarray*}
|I_{3}| &=& \sum_{k= - \infty}^{+\infty}  \chi(e^{-k} x^{w})\  w \int_{\frac{k}{w}}^{\frac{k+1}{w}} \big[ h \bigg(\frac{e^{u}}{x}\bigg) (u-\log(x))^{2} \big] \ du \\
&\leq& \sum_{|k-w \log(x)|< w \delta} \Big |\chi(e^{-k} x^{w})\  w \int_{\frac{k}{w}}^{\frac{k+1}{w}} \Big[ h \Big(\frac{e^{u}}{x}\Big) (u-\log(x))^{2} \Big] \ du \Big| \\&& + \sum_{|k-w\log(x)| \geq w \delta} \Big |\chi(e^{-k} x^{w})\  w \int_{\frac{k}{w}}^{\frac{k+1}{w}} \Big[ h \bigg(\frac{e^{u}}{x}\bigg) (u-\log(x))^{2} \Big] \ du \Big| \\
&:=& I_{3}^{'}+I_{3}^{''}.
\end{eqnarray*}
Since $\displaystyle \lim_{x \rightarrow 1} h(x)=0,$ we have $ \displaystyle |w I_{3}^{'}| \leq \frac{\epsilon}{3 w} (1+3 M_{2}(\chi)).$
Now using the fact that $h(x)$ is bounded, we obtain
\begin{eqnarray*}
|I_{3}^{''}| &\leq & \|h\|_{\infty} \sum_{|k-w\log(x)| \geq w \delta} \big |\chi(e^{-k} x^{w})\  w \int_{\frac{k}{w}}^{\frac{k+1}{w}} (u-\log(x))^2 \ du \big| \\
&\leq & \frac{(\epsilon+1)}{3 w^2} \|h\|_{\infty}.
\end{eqnarray*}
This gives $$ |w I^{''}| \leq \frac{(\epsilon+1)}{3 w} \|h\|_{\infty}.$$ Combining the estimates of $I_{1}-I_{3},$ we get the desired result.

\begin{Corollary}
The assumption that $f$ is bounded on $ \mathbb{R}^{+}$ can be relaxed by assuming that there are two positive constants $\alpha,\beta$ such that
$ |f(x)| \leq \alpha + \beta | \log (x)| , \hspace{0.5cm} \forall x \in \mathbb{R}^{+}.$
\end{Corollary}

\noindent\bf{Proof.}\rm \
First we show that the series (\ref{main}) is well defined for such $f.$ Indeed
\begin{eqnarray*}
|I_{w}^{\chi}f)(x)| &\leq &  \sum_{k= - \infty}^{+\infty} |\chi(e^{-k} x^{w})| \  w \int_{\frac{k}{w}}^{\frac{k+1}{w}} |f(e^{u}| \ du \\
& \leq & \sum_{k= - \infty}^{+\infty} |\chi(e^{-k} x^{w})| \  w \int_{\frac{k}{w}}^{\frac{k+1}{w}} (\alpha + \beta | u|) \ du \\
&\leq & (\alpha + \beta \log x) M_{0}(\chi)+ \beta  \sum_{k= - \infty}^{+\infty} |\chi(e^{-k} x^{w})| \  w \int_{\frac{k}{w}}^{\frac{k+1}{w}} |u - \log (x) | \ du \\
& \leq & \bigg(\alpha + \beta \log (x)+ \frac{\beta}{2w} \bigg) M_{0}(\chi)+ \frac{\beta}{w} M_{1}(\chi) \ < + \infty.
\end{eqnarray*}
This shows that the series $(I_{w}^{\chi}f)_{w>0}$ is absolutely convergent in $\mathbb{R}^{+}.$ Now for any fixed $a \in \mathbb{R}^{+},$ we define
$\displaystyle P_{2}(u):= f(a)+ (\theta f)(a) (u-\log(a)) + \frac{(\theta^{(2)} f)(a)}{2!} (u-\log(a))^{2}.$ From the Taylor's formula in terms of Mellin derivatives upto second order term, we can write as
$$ h \bigg(\frac{e^{u}}{a}\bigg) = \frac{f(e^{u}) - P_{2}(u)}{(u-\log (a))^{2}} \ ,$$
where $ h(.)$ is a bounded function such that $\displaystyle \lim_{t \rightarrow 1} h(t)=0.$ This implies that $h$ is bounded in the neighbourhood for $ |u - \log (a)| < \delta.$ For $ |u - \log (a)| \geq \delta,$ we have
\begin{eqnarray*}
|h(e^{u} a^{-1})| & \leq & \frac{|f(e^{u}|}{|u - \log (a)|^{2}} + \frac{|P_{2}(u)}{|u - \log (a)|^{2}} \\
& \leq & \frac{\alpha+ \beta |u|}{|u - \log (a)|^{2}} + \frac{|P_{2}(u)}{|u - \log (a)|^{2}} \ .
\end{eqnarray*}
This shows that $ h(.)$ is bounded on $ \mathbb{R}^{+}.$ Now we can proceed in the similar manner as in Theorem \ref{theorem2}, to get the same asymptotic formula.\\

The logarithmic modulus of continuity is defined by
$$ \omega(f,\delta):= \sup \{|f(x)-f(y)|:\  \mbox{whenever} \  |\log (x)-\log (y)| \leq \delta,\  \ \delta \in \mathbb{R}^{+}\} .$$
The properties of logarithmic modulus of continuity can be seen in \cite{bardaro9}. Now we obtain a quantitative estimate of the convergence of operator (\ref{main}) for $ f \in \mathcal{C}(\mathbb{R}^+).$

\begin{thm}\label{t3}
Let $ f \in \mathcal{C}(\mathbb{R}^+).$ Then, we have
 $$ |(I_{w}^{\chi}f)(x) - f(x)| \leq \lambda \ \omega \left( f,\frac{1}{w} \right),$$
where $\lambda = (M_{0}(\chi) +  M_{1}(\chi)).$
\end{thm}

\noindent\bf{Proof.}\rm \
We have
\begin{eqnarray*}
|(I_{w}^{\chi}f)(x) - f(x)| &\leq & \sum_{|k - w \log(x)|< w \delta}  |\chi(e^{-k} x^{w})| \  w \int_{\frac{k}{w}}^{\frac{k+1}{w}} |f(e^{u}) - f(x)| \ du    \\&& + \sum_{|k- w \log(x)| \geq w \delta} \big |\chi(e^{-k} x^{w})| \  w \int_{\frac{k}{w}}^{\frac{k+1}{w}} |f(e^{u}) - f(x)| \ du \\
&\leq& \sum_{|k-w \log(x)|< w \delta} \big |\chi(e^{-k} x^{w}) \big| \omega(f,|k - w \log(x)|) + 2 \|f\|_{\infty} M_{0}(\chi)\epsilon\\
&\leq&  \omega (f,\delta) \bigg( M_{0}(\chi)  + \frac{1}{w \delta} M_{1}(\chi) \bigg).
\end{eqnarray*}
Choosing $ \delta= \dfrac{1}{w}$ and $ \lambda= (M_{0}(\chi) +  M_{1}(\chi)),$ we get the desired estimate.

\section{Inverse result}
In this section, we derive an inverse result of approximation for the Kantorovich exponential sampling operators. In order to establish the saturation theorem, we first obtain a relation between $(S_{w}^{\chi}f)_{w>0}$ and $(I_{w}^{\chi}f)_{w>0}$ for $ f \in C^{(n)}(\mathbb{R}^{+}).$ This work is motivated from Kantorovich type generalized sampling operators studied in \cite{butzer2,costa4,costa5,bartoc}.

\begin{thm} \label{theorem4}
Let $f\in C^{(n)} (\mathbb{R^+}),$ $n\in \mathbb{N}.$ Then the following relation holds for every $ x \in \mathbb{R}^{+}:$
\begin{equation} \label{representation}
(I_w^\chi f)(x) = \sum_{j=0}^{n-1}\frac{( S_w^\chi \theta^{(j)}f))(x)}{(j+1)!} + \hat{R}_n^{w}(x),
\end{equation}
where $\displaystyle \hat{R}_n^{w}(x):= \frac{1}{n!} \sum_{k \in \mathbb{Z}} \chi(e^{-k}x^w) \ w \bigg[ \int_{k/w} ^{k+1/w}(\theta^{(n)} f)(\xi) (u-k/w)^n du \bigg] $ is absolutely convergent for $\xi \in (e^{k/w},e^{k+1/w})$  and $ w>0.$
\end{thm}

\noindent\bf{Proof.}\rm \
Using the $n^{th}$ order Mellin's Taylor's formula (\cite{butzer3}) and substituting $ x=e^{k/w} ,\  w>0$ in the formula and using the fact that $u \in \big[ k/w,(k+1)/w \big],$ we have $\xi \in (e^{k/w},e^{k+1/w}) .$
Thus we obtain
\begin{eqnarray*}
f(e^{u})&=& f(e^{k/w})+ (\theta f)(e^{k/w})(u-{k/w})+\frac{(\theta^{(2)} f)(e^{k/w})}{2!} (u-{k/w})^2 +...+\frac{(\theta^{(n-1)} f)(\xi)}{(n-1)!} (u-{k/w})^{n-1}\\ &&+R_n^{w},
 \, \, \, \mbox{where} \, \, \, \displaystyle R_n^{w} = \frac{(\theta^{(n)} f)(\xi)}{n!} (u-\log (x))^{n}, \ \ \ \ \xi \in (x,e^u).
\end{eqnarray*}
Now we evaluate \ $\displaystyle w\int_{k/w} ^{k+1/w} f(e^u) \ du $
\begin{eqnarray*}
&=& w\int_{k/w}^{k+1/w}\big[ f(e^{k/w})+ (\theta f)(e^{k/w})(u-{k/w})+...+R_n^{w}(u)\big]\\
&=& f(e^{k/w})+ (\theta f)(e^{k/w})\frac{1}{w (2!)}+...+w\int_{k/w}^{k+1/w}R_n^{w}(u) \ du.
\end{eqnarray*}
Using (\ref{main}), we obtain
\begin{eqnarray*}
(I_w^\chi f)(x)&=&(S_w^\chi f)(x) +\frac{1}{ (2!) w}(S_w^\chi (\theta f))(x)+...+\hat{R}_n^{w}(x)
=\sum_{j=0}^{n-1}\frac{( S_w^\chi \theta^{(j)}f)(x)}{(j+1)! \  w^j} + \hat{R}_n^{w}(x),
\end{eqnarray*}
where $\displaystyle\hat{R}_n^{w}(x):= \sum_{k\in \mathbb{Z}}\chi(e^{-k}x^w) w \int_{k/w} ^{k+1/w}\frac{(\theta^{(n)} f)(\xi)}{n!}(u-k/w)^n du. $ Next we show that the remainder term $\hat{R}_n^{w}$ is absolutely convergent in $ \mathbb{R}^{+}.$
\begin{eqnarray*}
\big| \hat{R}_n^{w} \big| &\leq & \frac{1}{{n!}} \sum_{k\in \mathbb{Z}}\big|\chi(e^{-k}x^w)\big| \bigg| w \int_{k/w} ^{k+1/w}(\theta^{(n)} f)(\xi) (u-k/w)^n du \bigg|
\leq \frac{\|\theta^{(n)}f\|_\infty}{{(n+1)!}} \frac{1}{w^n} M_0(\chi)< +\infty .
\end{eqnarray*}
This completes the proof.\\

Now, we have the following inverse result for the Kantorovich exponential sampling operators $(I_{w}^{\chi}f)_{w>0}.$

\begin{thm} \label{t5}
Let $\chi$ be the kernel function such that $m_{1}(\chi,u)=0.$ Suppose that $ f \in C^{(2)}(\mathbb{R}^{+})$ and
$$ \| I_{w}^{\chi}f - f \|_{\infty} = o(w^{-1}) \hspace{0.5 cm} \mbox{as} \ \ w \rightarrow \infty.$$ Then $f$ is  constant on $\mathbb{R}^{+}.$
\end{thm}

\noindent\bf{Proof.}\rm \
Using the relation (\ref{representation}) for $n=1,$ we have
$$ |(I_{w}^{\chi}f)(x) - f(x)|= |(S_{w}^{\chi}f)(x)- f(x) + \hat{R}_{1}^{w}(x)|.$$
Using the assumption that $ \| I_{w}^{\chi}f - f \|_{\infty} = o(w^{-1}),$ we obtain
$$ |(S_{w}^{\chi}f)(x)- f(x) + \hat{R}_{1}^{w}(x)| = o(w^{-1}),$$
which implies that
\begin{equation} \label{inverse}
\lim_{w\rightarrow \infty} w [ S_{w}^{\chi}f)(x)- f(x)] + \lim_{w\rightarrow \infty} w \ \hat{R}_{1}^{w}(x)=0.
\end{equation}
Consider
\begin{eqnarray*}
| 2 w \ \hat{R}_{1}^{w}(x) - (\theta f)(x)| &=& |2 w \ \hat{R}_{1}^{w}(x) - (S_{w}^{\chi}\theta f)(x)+ (S_{w}^{\chi}\theta f)(x) - (\theta f)(x) | \\
&\leq & |(S_{w}^{\chi}\theta f)(x) - (\theta f)(x) | + |2 w \ \hat{R}_{1}^{w}(x)-(\theta f)(x)|:=I_{1}+I_{2}.
\end{eqnarray*}
Using Theorem 5 and Corollary 2 in \cite{bardaro7}, we have $|I_{1}| \leq \epsilon.$ Now we estimate $I_{2}.$
\begin{eqnarray*}
| I_{2}|
& \leq & \bigg| \bigg[ 2w \  \sum_{k= - \infty}^{+\infty}  \chi(e^{-k} x^{w}) \  w \int_{\frac{k}{w}}^{\frac{k+1}{w}} (\theta f)(\xi) \left(u- \frac{k}{w} \right) \ du \bigg] - (S_{w}^{\chi} \theta f)(x) \bigg|\\
&\leq &  2w^{2} \bigg| \bigg[ \  \sum_{k= - \infty}^{+\infty} \chi(e^{-k} x^{w})\   \int_{\frac{k}{w}}^{\frac{k+1}{w}}  \left(u- \frac{k}{w} \right) \Big[ (\theta f)(\xi) - (\theta f) (e^{\frac{k}{w}}) \Big] \ du \bigg] \bigg| .
\end{eqnarray*}
Since $\xi \in \big[\frac{k}{w}, \frac{k+1}{w}]$ and by using the continuity of $ (\theta f)(x),$ we have
$$ |(\theta f)(\xi) - (\theta f) (e^{\frac{k}{w}})| < \epsilon$$
which gives $ |I_{2}| < \epsilon M_{0}(\chi).$ Combining the estimates $ I_{1}-I_{2},$ we obtain
$$ \lim_{w\rightarrow \infty} [ 2w \hat{R}_{1}^{w}(x)-(\theta f)(x)] =0.$$
Using the above estimate and (\ref{inverse}), we have $ (\theta f)(x)=0$  $\forall x \in \mathbb{R}^{+}.$ This  implies that $f$ is constant on $\mathbb{R}^{+}.$

\section{Examples of the kernels and Graphical representation}
In this section, we present few examples of the kernel functions based on the theory of Mellin's transform which also satisfy the assumptions of the presented theory. The first example of kernel functions in this direction is the family of Mellin- B spline kernels \cite{bardaro7}.

\subsection{B-splines in the Mellin setting}
The B-splines of order $n$ in the Mellin setting for $ x \in \mathbb{R}^{+}$ are defined as
$$\bar{B}_{n}(x):= \frac{1}{(n-1)!} \sum_{j=0}^{n} (-1)^{j} {n \choose j} \bigg( \frac{n}{2}+\log(x)-j \bigg)_{+}^{n+1}.$$
$\bar{B}_{n}(x) $ is compactly supported for every $ n \in \mathbb{N}.$ The Mellin transformation of $\bar{B}_{n}$ (see \cite{bardaro7}) is given by
\begin{eqnarray} \label{splinefourier}
\hat{M}[\bar{B}_{n}](c+iw) = \bigg( \frac{\sin(\frac{w}{2})}{(\frac{w}{2})} \Bigg)^{n},  \ \ \hspace{0.5cm} w \neq 0.
\end{eqnarray}
To show that the above kernel satisfies the assumptions (i)-(ii), we use the Mellin's - Poisson summation formula \cite{butzer4}) which has the form
$$ (i)^{j} \sum_{k= - \infty}^{+\infty} \chi(e^{k} x) ( k-\log(u))^{j} = \sum_{k= - \infty}^{+\infty} \hat{M}^{(j)}[\chi](2k \pi i) \ x^{-2 k \pi i}, \ \ \ \ \ \ \mbox{for} \ k \in \mathbb{Z}.$$
The following lemma will be useful in this direction.
\begin{lem}\label{l1}
The condition $\displaystyle \sum_{k= - \infty}^{+\infty}  \chi(e^{-k} x^{w})=1$ is equivalent to the condition
\begin{equation} \label{e4}
\hat{M}[\chi] (2k \pi i) =
     \begin{cases}
     {1,} &\quad\text{if,} \ \ \ \ {k=0} \\
     {0,} &\quad\text{if,} \ \ \ \ {\text{otherwise}}\\
   \end{cases}
\end{equation}

Moreover, the condition $ m_{j}(\chi,u)=0$ for $ j= 1,2,...n$ is equivalent to the condition \\ $ \hat{M^{(j)}}[\chi](2k \pi i)= 0$ for $ j= 1,2...,n$ and $ \forall \ k \in \mathbb{Z}.$
\end{lem}
\noindent\bf{Proof.}\rm \
The proof can be obtained by the similar arguments as given in the proof of Lemma 2 and Lemma 3 in \cite{butzer2}.

Again from (\ref{splinefourier}) we have,
\begin{equation*}
\hat{M}[\bar{B}_{n}] (2k \pi i) =
     \begin{cases}
      {1,} &\quad\text{if,} \ \ \ \ {k=0} \\
     {0,} &\quad\text{if,} \ \ \ \ {\text{otherwise.}}\\
   \end{cases}
\end{equation*}
This shows that $ \bar{B}_{n}(x)$ satisfies the condition (i). Since $ \bar{B}_{n}(x)$ is compactly supported, the condition (ii) is also satisfied. By using Lemma 1 for $j=1,$ we obtain that the first order moment for $ \hat{M}[\bar{B}_{n}]$ vanishes in $\mathbb{R}^{+}$ and hence it satisfies the assumption of Theorem \ref{theorem2}. 
First we show the approximation of $ f_{1}(x)$ by $(I_{w}^{\chi}f_{1})_{w>0},$ where
\begin{equation*}
f_{1}(x) =
     \begin{cases}

     {0 ,} &\quad\text{} \ \ \ \ { \frac{1}{2} \leq x <1}\\
       {\frac{-2}{x},} &\quad\text{} \ \ \ \ { 1 \leq x <4. }\\
   \end{cases}
\end{equation*}

\begin{figure}[h]
\centering
{\includegraphics[width=0.8\textwidth]{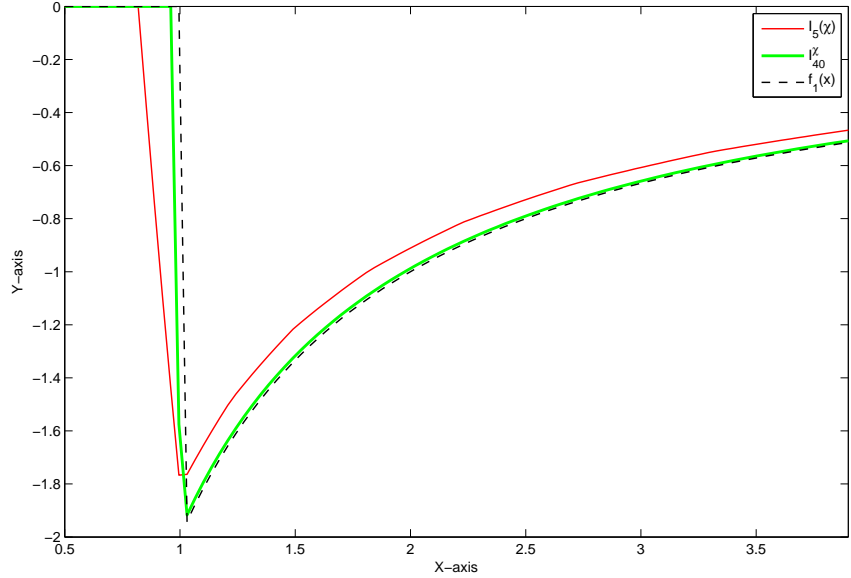}}
\caption{This figure exhibits the approximation of $f_{1}(x)$ (Black) by the series $(I_{w}^{\chi}f_{1})_{w>0}$ for $ w= 5, 40$ (Red and Green respectively) based on Mellin B-spline kernel.}
\end{figure}

\begin{tb}\label{table1}\centering
 {\it Error estimation (upto $4$ decimal points) in the approximation of $f_{1}(x)$ by $I_{w}^{\chi}f_{1}(x)$ for $w=5,40,70.$}

$  $

\begin{tabular}{|l|l|l|l|}\hline
 $x$&$ |f_{1}(x) - I_{5}^{\chi}f_{1}(x)|$&$|f_{1}(x)-I_{40}^{\chi}f_{1}(x)|$&$|f_{1}(x)-I_{70}^{\chi}f_{1}(x)|$\\
 \hline
 $1.1$&$0.1621$ & $0.0225$ & $0.0129$\\
  \hline
 $1.8$ & $0.1028$  &  $0.0137$ & $0.0079$ \\
  \hline
 $2.9$ & $0.0620$  &  $0.0085$ & $0.0049$\\
 \hline
  $3.8$ & $0.0471$  &  $0.0065$ & $0.0037$  \\  \hline

             \end{tabular}
   \end{tb}

\newpage

\subsection{Mellin's-Fejer kernel}
The general form of Mellin's-Fejer kernel is given by
$$ F_{\alpha}^{c}(x):= \frac{\alpha}{2 \pi x^
{c}} sinc^{2} \bigg(\frac{\alpha}{\pi} \log(\sqrt{x}\bigg) ,$$
where $ c \in \mathbb{R}, \alpha >0$ and $ x \in \mathbb{R}^{+}.$ The Mellin's transform for $ F_{\alpha}^{c}$ is given by
\begin{equation*}
\hat{M}[F_{\alpha}^{c}] (c+i w) =
     \begin{cases}
     \text      {1 - $\frac{|w|}{\alpha}$,} &\quad\text{if,} \ \ \ \ {|w| \leq \alpha} \\
     \text{0,} &\quad\text{if,} \ \ \ \ {|w|>\alpha}.\\
   \end{cases}
\end{equation*}
The Mellin-Fejer kernel also satisfies the assumption (i)-(ii) analogously (see \cite{bardaro7}), but it fails to satisfy the moment condition of Theorem 3, i.e, $ m_{1}(\chi,u)=0.$ This can be seen by using the equivalent condition mentioned in Lemma \ref{l1} for $ j=1.$
Next we present the approximation of the function $f_{2}(x)$ by the operator $(I_{w}^{\chi}f_{2})_{w>0},$ where
\begin{equation*}
f_{2}(x) =
     \begin{cases}

     {\cos x,} &\quad\text{} \ \ \ \ { 1 \leq x < 4}\\
       {0,} &\quad\text{} \ \ \ \ { 0 \leq x <1.} \\
   \end{cases}
\end{equation*}

\begin{figure}[h]
\centering
{\includegraphics[width=1.0 \textwidth]{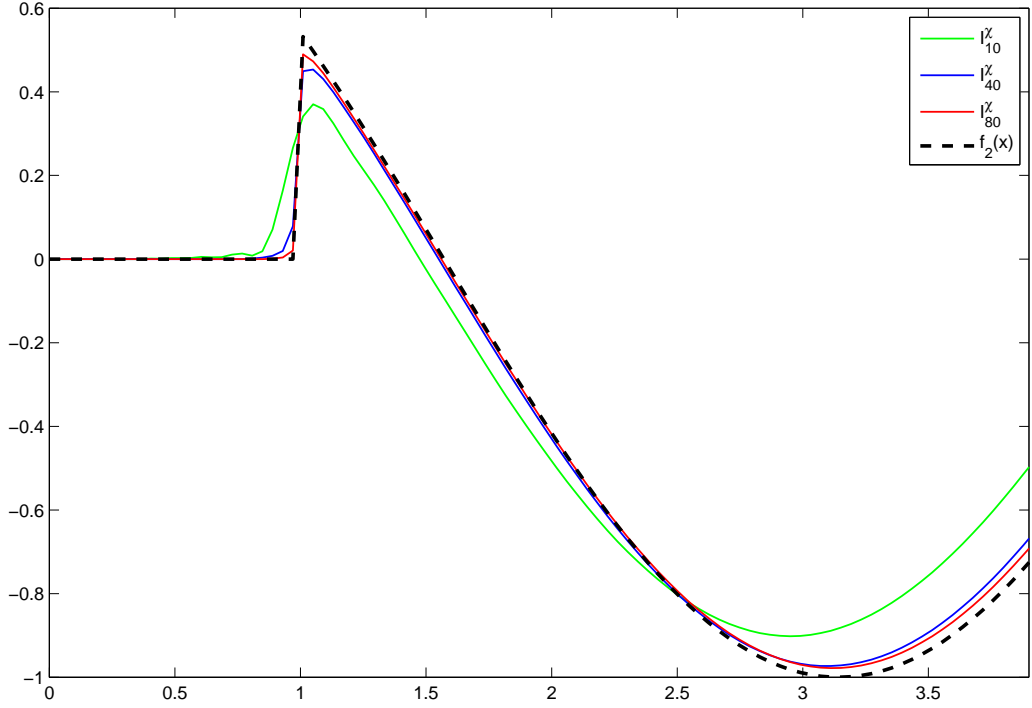}}
\caption{This figure exhibits the approximation of $f_{2}(x)$ (Black) by the series $(I_{w}^{\chi} f_{2})_{w>0}$ for $ w= 10, 40,80$ (Green,Blue and Red respectively) based on Mellin-Fejer kernel $F_{\pi}^{0}(x).$}
\end{figure}
\newpage
\begin{tb}\label{table2}\centering
 {\it Error estimation (upto $4$ decimal points) in the approximation of $f_{2}(x)$ by $I_{w}^{\chi}f_{2}(x)$ for $w=10,40,80.$}

$  $

\begin{tabular}{|l|l|l|l|}\hline
 $x$&$ |f_{2}(x) - I_{10}^{\chi}f_{1}(x)|$&$|f_{2}(x)-I_{40}^{\chi}f_{2}(x)|$&$|f_{2}(x)-I_{80}^{\chi}f_{2}(x)|$\\
 \hline
 $1.4$&$0.0954$ & $0.0216$ & $0.0123$\\
  \hline
 $2.3$ & $0.0322$  &  $0.0059$ & $0.0033$ \\
  \hline
 $3.4$ & $0.1635$  &  $0.0391$ & $0.0271$\\
 \hline
  $3.9$ & $0.2262$  &  $0.0571$ & $0.0336$  \\  \hline

             \end{tabular}
   \end{tb}

\end{document}